\documentclass[12pt]{amsart} 
\usepackage{amssymb}

\DeclareMathAlphabet{\mat}{U}{msb}{m}{n}
\DeclareMathAlphabet{\cur}{U}{eus}{m}{n}
\DeclareMathAlphabet{\got}{U}{euf}{m}{n}

\newcommand{\lra}{\longrightarrow}
\newcommand{\dis}{\displaystyle}

\DeclareMathOperator{\Ad}{Ad}
\DeclareMathOperator{\ad}{ad}

\begin{document}
\title[Sur la structure transverse...]{SUR LA STRUCTURE TRANSVERSE A UNE ORBITE NILPOTENTE ADJOINTE.}
\author{Herv\'e Sabourin }

\address
{  UMR 6086 (GLG) CNRS, D\'epartement de Math\'ematiques,
T\'el\'eport 2, BP 30179,
Boulevard Marie et Pierre Curie,
86962 Futuroscope-Chasseneuil Cedex, France}

\email{ sabourin@mathlabo.univ-poitiers.fr}
 
\begin{abstract}We are interested in  Poisson structures transverse
to nilpotent adjoint
orbits in a complex semi-simple Lie algebra, and we
study their
polynomial nature, introduced in ([CU-RO]). Furthermore, in the case
of $sl_n$, 
we construct some families of nilpotent orbits with quadratic
transverse structures.   
\end{abstract}

\maketitle

\section*{\bf 0.Introduction}
L'objet de ce travail est une contribution \`a l'\'etude de la structure
de Poisson transverse \`a une orbite nilpotente adjointe dans une
alg\`ebre de Lie semi-simple complexe.

Soit $\got g$ une telle alg\`ebre de Lie, $G$ son groupe adjoint. Chaque
$G$-orbite nilpotente $O$ est enti\`erement caract\'eris\'ee par la donn\'ee d'un
$sl_{2}$-triplet $(h,e,f)$, o\`u $e$ est un g\'en\'erateur de l'orbite et
$h$ sa caract\'eristique. L'alg\`ebre $\got g$ est identifi\'ee \`a son dual
${\got g}^{*}$, via la forme de Killing, et peut ainsi \^etre munie de
la structure de Lie-Poisson correspondante. On sait alors r\'ealiser la
structure de Poisson transverse \`a $O$ dans $\got g$, en choisissant un
suppl\'ementaire quelconque $\got n$ du centralisateur $ {\got g}^e$ de
$e$ dans $\got g$ et en consid\'erant la vari\'et\'e affine $N = e + {\got
n}^{\perp}$. On sait, d'apr\`es A.Weinstein ([WE]), que les structures
de Poisson transverses associ\'ees \`a deux suppl\'ementaires sont
n\'ecessairement isomorphes. 

 Dans ([SG]), M.Saint-Germain a remarqu\'e que si l'on d\'efinit un syst\`eme
de coordonn\'ees locales $(q) = (q_{1},\dots,q_{n})$ pour $N$, alors la
structure de Poisson est \`a coefficients rationnels en la variable
$(q)$. Dans ([CU-RO]), R.Cushman et M.Roberts ont consid\'er\'e le
suppl\'ementaire $ {\got n} = Im \ad f$ et ont montr\'e que, pour un tel
choix, la structure de Poisson transverse \'etait polyn\^omiale en la
variable $(q)$. 

Nous consid\'erons alors l'ensemble $ {\mathcal N}_{h}$ des
suppl\'ementaires $\ad h$-invariants de $ {\got g}^e$ et montrons que,
pour tout \'el\'ement de cet ensemble, la structure transverse
correspondante est encore polyn\^omiale. 

Nous proposons, par ailleurs, dans l'exemple de l'orbite sous-r\'eguli\`ere de
$sl(4,{\mathbb C})$, un suppl\'ementaire qui n'est pas $h$-invariant
pour lequel la structure transverse n'est pas polyn\^omiale, ce
qui justifie l'introduction de la vari\'et\'e $ {\mathcal N}_{h}$. 

A chaque suppl\'ementaire $h$-invariant $\got n$, on peut donc
associer le degr\'e $deg_{N}$ de la structure transverse correspondante,
d\'efini comme \'etant le 
le degr\'e maximum des polyn\^omes intervenant dans la
structure de Poisson. Naturellement, la question se pose de savoir si
cette notion de degr\'e peut \^etre d\'efinie de mani\`ere intrins\`eque, c'est-\`a-dire
ind\'ependamment du choix du suppl\'ementaire choisi
dans $ {\mathcal N}_{h}$. 

Nous consid\'erons \`a nouveau, \`a ce sujet, l'exemple de l'orbite sous-r\'eguli\`ere de
$sl(4,{\mathbb C})$ et donnons deux
suppl\'ementaires $h$-invariants dont les structures transverses sont de 
``degr\'e'' distincts, ce qui r\'epond n\'egativement \`a la question pos\'ee.

Nous nous int\'eressons enfin au cas des structures transverses
``quadratiques'', c'est-\`a-dire des structures polyn\^omiales de
degr\'e au
plus \'egal \`a $2$. Selon la terminologie introduite par M.Roberts
([RO]), on consid\`ere les orbites conormales, celles pour lesquelles le
centralisateur $ {\got g}^e$ poss\`ede un suppl\'ementaire qui est une
sous-alg\`ebre. La structure transverse correspondant \`a un tel
suppl\'ementaire est alors quadratique ([OH]). Nous montrons dans le
paragraphe 3 que toute orbite nilpotente ``sph\'erique'' est une orbite
conormale et, dans le cas $sl_n({\mathbb C})$, nous construisons une
famille d'orbites nilpotentes conormales, contenant entre autres les
orbites sph\'eriques et l'orbite r\'eguli\`ere, g\'en\'eralisant en cela
certains r\'esultats de M.Rais ([RA]).

L'exemple de l'orbite sous-r\'eguli\`ere de $sl_4({\mathbb C})$ nous
permet d'exhiber une structure transverse quadratique associ\'ee \`a un
suppl\'ementaire qui n'est pas une sous-alg\`ebre. R\'eciproquement, on peut
se demander si toute structure transverse d'une orbite conormale est
quadratique; nous donnons pour finir un exemple dans $sl_5({\mathbb
C})$ qui r\'epond n\'egativement \`a cette question.  

Je tiens \`a remercier tout particuli\`erement M.Rais et M.Duflo pour les
fr\'equents  \'echanges ou nombreuses 
discussions fort utiles que nous avons eus sur ce sujet. Je tiens aussi
\`a remercier M.Saint-Germain pour l'attention qu'il a su
porter \`a ce travail.

\section{\bf Structure transverse \`a une orbite coadjointe.}
\subsection{}
Soit $\got g$ une alg\`ebre de Lie complexe, $G$ son groupe
adjoint, ${\got g}^{*}$ son dual. On sait que ${\got g}^{*}$ peut \^etre
munie d'une structure de 
vari\'et\'e de Poisson et on notera $\Lambda$ le tenseur de Lie-Poisson
correspondant . En tout point $\mu \in {\got g}^{*}$, on identifie l'espace
cotangent $T^{*}_{\mu}{\got g}^{*}$ avec $\got g$  et, moyennant cette
identification, on d\'efinit $\Lambda$ par :
$$
\forall X,Y \in {\got g}, \Lambda_{\mu}(X,Y) = \mu([X,Y])
$$
 
Les feuilles symplectiques dans ${\got g}^{*}$ sont les
orbites coadjointes et, selon A.Weinstein ([WE]), on peut donc d\'efinir
pour chaque orbite $G.\mu$ une
 structure transverse, c'est-\`a-dire une sous-vari\'et\'e $N$ de ${\got g}^{*}$
contenant $\mu$,
telle qu'en chaque point $\nu$ d'un voisinage $U$ de $\mu$ dans $N$, on ait:
$$
T_\nu{\got g}^{*} = T_\nu(G.\mu) + T_{\nu}N
$$
Selon A.Weinstein, on sait de plus que cette structure transverse est
unique \`a un isomorphisme de vari\'et\'es de Poisson pr\`es.

\subsection{}
  
Pour construire une telle structure transverse, on proc\`ede de la
mani\`ere suivante (voir M.Saint-Germain ([SG]) : On consid\`ere le centralisateur $ {\got g}^{\mu}$ de
$\mu$ dans $\got g$ et un suppl\'ementaire quelconque $\got n$ de $ {\got
g}^{\mu}$. On note $ {\got n}^{0}$ (respectivement ${\got g}^{\mu,0}
$)l'annulateur de $\got n $ (respectivement de $ {\got g}^{\mu}$) dans $
{\got g}^{*}$ et on a :
$$
{\got g}^{*} = {\got g}^{\mu,0} \oplus {\got n}^{0}
$$
 
On pose : $N = \mu + {\got n}^{0}$.

On voit facilement que :
$$
T_{\mu}(G.\mu) = {\got g}.\mu = {\got g}^{\mu,0}, \ T_{\mu}N = {\got n}^{0}
$$
Ceci permet de constater que $N$ est bien transverse \`a $G.\mu$ au
point $\mu$ et on v\'erifie que cette propri\'et\'e reste vraie dans un
voisinage bien choisi de $\mu$ dans $N$.

Pour d\'eterminer le tenseur de Poisson $\Lambda_{N}$ de la
vari\'et\'e de
Poisson $N$ ainsi d\'efinie, on peut utiliser la formule dite ``formule
des contraintes de P.A.M. Dirac'', donn\'ee par P.G.Bergmann et
I.Goldberg ([BE-GO]). 

Soit $(Z_{i})_{1\leq i\leq k}$ une base de $ {\got g}^{\mu}$ et
$(X_{j})_{1\leq j\leq p}$ une base du suppl\'ementaire $\got n$. Un
syst\`eme de coordonn\'ees lin\'eaires pour un \'el\'ement quelconque $\varphi$
de ${\got g}^{*} $ est alors donn\'e par :
$$
z_{i}(\varphi) = \varphi(Z_{i}), x_{j}(\varphi) = \varphi(X_{j}),
\forall (i,j), 1 \leq i \leq k, 1 \leq j \leq p
$$
En particulier, on notera : $\mu = (\mu_{1}, \dots, \mu_{n})$ les
coordonn\'ees de $\mu$.

Soit $U$ un ouvert de ${\got g}^{*}$. On a :
$$
N \cap U = \{\varphi \in U, x_{j}(\varphi) = \mu_{j}, \forall j, 1
\leq j \leq p \}
$$
 
Ainsi, l'application $\varphi \lra (z_{1}(\varphi),\dots,z_{k}(\varphi))$ est
une carte de $N$.

Soit $(i,j,l,m), 1 \leq i,j \leq k, 1 \leq l,m \leq p$ et $\nu \in {\got
n}^{0}$. 
On pose :
\begin{align*}
C_{N}(\nu)_{l,m} &= (\mu + \nu)([X_{l},X_{m}]) \\
D_{N}(\nu)_{l,j} &= \nu([X_{l},Z_{j}])   \\
A_{N}(\nu)_{i,j} &= \nu([Z_{i},Z_{j}]) 
\end{align*}
On notera respectivement $C_{N}(\nu), D_{N}(\nu)$ et $A_{N}(\nu)$ les
matrices de termes g\'en\'eraux donn\'es par les formules pr\'ec\'edentes. 
On sait que la matrice $C_{N}(\nu)$ est inversible lorsque $\nu$
parcourt un voisinage $V$ de $0$, bien choisi dans ${\got n}^{0} $. 

On pose, de plus : $B_{N}(\nu) = \ ^{t}D_{N}(\nu)$. 

La matrice du tenseur
$\Lambda_{N}$ en un point quelconque $(\mu + \nu)$ de $\mu + V$ sera not\'ee
$\Lambda_{N}(\nu)$. Son terme g\'en\'eral est le crochet de Poisson
$\{z_{i},z_{j}\}(\mu + \nu)$.

La formule des contraintes de Dirac est alors la suivante :
\begin{equation}
\Lambda_{N}(\nu) = A_{N}(\nu) + B_{N}(\nu)C_{N}(\nu)^{-1}D_{N}(\nu)
\end{equation}

On constate finalement que les coefficients $D_{N}(\nu)_{l,j} $  
et $ A_{N}(\nu)_{i,j}$ sont lin\'eaires en les coordonn\'ees $(z_{j}(\nu))$,
que le coefficient $ C_{N}(\nu)_{l,m}$ est affine en les $(z_{j}(\nu))$. Il
s'en suit que la matrice inverse de $C_{N}(\nu)$ a des coefficients
rationnels en les coordonn\'ees $(z_{j}(\nu))$. Ainsi, chaque crochet de
Poisson $\{z_{i}, z_{j}\}$ est rationnel en $(z_{1},\dots,z_{k})$.

\section{\bf Le cas d'une alg\`ebre de Lie semi-simple.}

\subsection{}
On suppose maintenant que $\got g$ est une alg\`ebre de Lie semi-simple
complexe, que $\got h$ est une sous-alg\`ebre de Cartan de $\got g$. On
d\'efinit un syst\`eme de racines $\Delta({\got h})$, $\Pi({\got h})$ une
base de racines simples  et on note $K$ la forme de Killing qui permet d'identifier
$\got g$ \`a son dual $ {\got g}^{*}$.
Soit $x \in {\got g}$ et $\got n$ un suppl\'ementaire quelconque de $
{\got g}^{x}$ dans $\got g$. On note $ {\got n}^{\perp}$
l'orthogonal dans $ {\got g}$ de $\got n$ relativement \`a $K$. Alors, $N = x + {\got n}^{\perp}$
est la structure transverse \`a l'orbite adjointe $G.x$.  

On peut d\'efinir comme pr\'ec\'edemment le tenseur de Poisson $\Lambda_{N}$
\`a partir des bases $(Z_{i}), (X_{j})$. 
On consid\`ere la base $(\overline{Z_{i}}, \overline {X_{j}})$, duale de
la base pr\'ec\'edente relativement \`a $K$. N\'ecessairement, la famille
$(\overline{Z_{i}})$ est une base de $ {\got
n}^{\perp}$ et on identifie $ {\got n}^{\perp}$ \`a $ {\mathbb C}^{k}$ par l'application $\dis{q = (q_{1},\dots,q_{k}) \lra
\sum_{i=1}^{i=k}q_{i}\overline{Z_{i}}} $. 

On reprend ensuite les notations pr\'ec\'edentes, soit :
\begin{align*}
C_{N}(q)_{l,m} &= \sum_{s=1}^{s=k}K(x + q_{s}\overline{Z_{s}},[X_{l},X_{m}]) \\
D_{N}(q)_{l,j} &=
\sum_{s=1}^{s=k}K(q_{s}\overline{Z_{s}},[X_{l},Z_{j}]) \\ 
A_{N}(q)_{i,j} &=
\sum_{s=1}^{s=k}K(q_{s}\overline{Z_{s}},[Z_{i},Z_{j}]) 
\end{align*}
La matrice $\Lambda_{N}(q)$ du tenseur $\Lambda_{N}$ est encore
d\'efinie par la formule (1) et ses coefficients sont des fonctions
rationnelles en $(q_{1}, \dots,q_{k})$.

\subsection{}
Nous allons rappeler, dans cette section,  quelques faits et
notations sur les orbites nilpotentes utiles pour la suite. On pourra,
pour plus de d\'etails sur ce sujet, se r\'ef\'erer \`a ([SP-ST]).

A chaque $G$-orbite nilpotente $O$, on sait associer sa classe de
conjugaison de Jacobson-Morosov de $sl_{2}$-triplets. Dans cette
classe, il existe un unique 
$sl_{2}$-triplet $(h^{+},e,f)$  tel que :
 
$$
\forall \alpha_{i} \in \Pi({\got
h}), \alpha_{i}(h^{+}) \in \{0,1,2\}
$$
La suite
$(\alpha_{i}(h^{+}))$ est la {\it caract\'eristique} de l'orbite. Dans
ces conditions, l'\'el\'ement $e$ est un g\'en\'erateur de $O$.

on consid\`ere alors la d\'ecomposition de $\got g$ en sous-espaces
propres suivant l'action de $\ad h^{+}$, soit :
$$
{\got g} = \bigoplus_{m \in {\mathbb Z}}{\got g}(m)
$$
Le plus grand entier strictement positif $i$ tel que $ {\got g}(i)
\not= 0$ est la {\it hauteur} de l'orbite et sera not\'e $h(e)$.

Pour les alg\`ebres de Lie simples classiques de type $sl(V), so(V)$ ou $ sp(V)$, il est souvent commode de
d\'ecrire les orbites nilpotentes en termes de partition de l'entier $n
= \dim V$. En particulier, l'ensemble des orbites nilpotentes de
$sl_{n}({\mathbb C})$ s'identifie \`a l'ensemble des partitions
$(p_{1},\dots,p_{s}) , p_{1} \geq p_{2} \geq \dots \geq p_{s}$ avec
$\dis{\sum_i p_{i} = n}$. La hauteur d'une telle orbite est \'egale \`a
$2(p_{1}-1)$.     

Enfin, suivant des notations usuelles, \`a chaque racine $\alpha \in
\Delta({\got h})$,
on fait correspondre les vecteurs-racines
$(H_{\alpha},X_{\alpha},X_{{-\alpha}}) $ d'une base de Chevalley de $
\got g$.

\subsection{}
On consid\`ere donc le cas d'une orbite adjointe nilpotente $O = G.e$
et son $sl_{2}$-triplet associ\'e
$(h^{+},e,f)$, avec $h^{+} \in {\got h}$. Soit $ {\got s}$ la sous-alg\`ebre de Lie de $ \got
g$, engendr\'ee par ce $sl_{2}$-triplet. 

Outre la d\'ecomposition de $ \got g$ en sous-espaces propres , sous l'action de $\ad h^{+}$,
 on a aussi une d\'ecomposition de $\got g$ en composantes irr\'eductibles
 sous l'action de $\got s$, soit :
$$
{\got g} = \bigoplus_{i = -h(e)}^{i = h(e)} {\got g}(i) = \bigoplus_{k} E_{n_{k}}
$$
Chaque $ {\got g}(i)$ est l'espace propre associ\'e au poids entier  $i$, chaque
$E_{n_{k}}$ est la composante irr\'eductible, suivant l'action de $\got
s$ dans $ \got g$, de plus haut poids $n_{k}$.
Le centralisateur $ {\got g}^{e}$ est alors engendr\'e par les vecteurs
de plus haut poids de cette d\'ecomposition.
  
On note $ {\mathcal N}_{h}$ l'ensemble des suppl\'ementaires $\ad
h^{+}$-invariants de $ {\got g}^{e}$ dans $\got g$. Pour chaque $ \got n$
dans l'ensemble $ {\mathcal N}_{h}$, on a la d\'ecomposition :
$$
{\got n} = \bigoplus_{i \in {\mathbb Z}} {\got n}(i), \  {\got n}(i) =
{\got g}(i) \cap {\got n}
$$
On va voir, dans ce qui va suivre, le r\^ole important jou\'e par l'ensemble
$ {\mathcal N}_{h}$ dans l'\'etude de la structure transverse \`a l'orbite
$O$. Un premier r\'esultat est relatif \`a la structure de vari\'et\'e
de $
{\mathcal N}_{h}$.

{\bf Proposition 2.1} {\it L'espace $ {\mathcal N}_{h}$ est une
vari\'et\'e 
alg\'ebrique irr\'eductible de dimension :
$$
\dis{2 \sum_{i=0}^{i=h(e)}
\dim {\got g}^e(i)(\dim {\got g}(i) - \dim {\got g}^e(i))}
$$}

On rappelle tout d'abord le r\'esultat classique suivant, ainsi que
sa d\'emonstration (voir [DI], ch.1.11):

{\bf Proposition 2.2 :} 
{\it Soit $V $un $ {\mathbb C}$-espace vectoriel de dimension $n$, $W$
un sous-espace de $V$ de dimension $d$. Alors, 
l'ensemble $\mathcal S$ des suppl\'ementaires de $W$ dans $V$ est un
ouvert affine de $Gr(V,n-d)$, la
grassmanienne de dimension $n-d$ dans $V$.} 

Il est bien connu, en effet, que $ {\mathcal S}$ est un ouvert de
$Gr(V,n-d)$.

Consid\'erons une base $(e_{1},\dots,e_{n})$ de $V$ telle que
$(e_{1},\dots,e_{d})$ soit une base de $W$. Soit $F \in {\mathcal
S}$. Alors, pour tout $i, d+1 \leq i \leq n$, il existe des scalaires
$(a_{i,j}, 1 \leq j \leq d)$ tels que : $\dis{e_{i} + \sum_{j} a_{ij}e_{j}
\in F}$. Il s'en suit que l'application $F \lra (a_{ij})_{\substack{d+1
\leq i \leq n \\ 1 \leq j \leq d}}$ est un isomorphisme de vari\'et\'es de
$\mathcal S$ sur $ {\mathbb C}^{d(n-d)}$. Ainsi, $\mathcal S$ est un
ouvert affine de $Gr(V,n-d)$ et donc une 
vari\'et\'e alg\'ebrique irr\'eductible de dimension $d(n-d)$.

{\bf Preuve de la proposition 2.1 :} 
L'ensemble $ {\mathcal N}$ des suppl\'ementaires de $ {\got
g}^e$ dans $\got g$ est une vari\'et\'e alg\'ebrique irr\'eductible dont $
{\mathcal N}_{h}$ est une sous-vari\'et\'e.

Soit $i$ un entier, $-h(e) \leq i \leq h(e)$. D\'esignons par $
{\mathcal N}_{i}$ la vari\'et\'e des suppl\'ementaires de $ {\got g}^e(i) =
{\got g}^e \cap {\got g}(i)$ dans $ {\got g}(i)$. En utilisant ce qui
pr\'ec\`ede, on peut affirmer que $ {\mathcal N}_{i}$ est une
vari\'et\'e 
alg\'ebrique irr\'eductible de dimension :  $\dim {\got g}^e(i)(\dim {\got
g}(i) - \dim {\got g}^e(i))$. L'application $\dis{ {\got n} \lra
\bigoplus_{i=-h(e)}^{i=h(e)} {\got n}(i)}$ est alors un isomorphisme de
vari\'et\'es de $ {\mathcal N}_{h}$ sur $\dis{\prod_{i=-h(e)}^{i=h(e)}
{\mathcal N}_{i}}$. Ceci permet de d\'emontrer le r\'esultat souhait\'e.

{\bf Th\'eor\`eme 2.3 :} {\it Soit $ {\got n} \in {\mathcal
N}_{h}$. Alors, tous les coefficients de la matrice $\Lambda_{N}(q)$
sont des fonctions polyn\^omiales en la variable $q$.}

{\bf Preuve :}
ce r\'esultat a \'et\'e d\'emontr\'e par R.Cushman et M.Roberts
dans ([CU-RO]), dans le cas o\`u $ {\got n} = $Im ad $f$. Les arguments
employ\'es peuvent s'adapter sans difficult\'es \`a cette situation plus
g\'en\'erale.

On peut ausi utiliser les arguments donn\'es par P.Slodowy ([SL]). En
effet, il suffit 
de v\'erifier que le d\'eterminant de la matrice $C_{N}(q)$ est
ind\'ependant de la variable $q$. Pour cela, on utilise une
action $\rho$ de ${\mathbb C}^{*} $ sur la vari\'et\'e $N$, via
l'\'el\'ement $h^+$, introduite par T.A.Springer et R.Steinberg dans
([SP-ST], paragraphe 4) et d\'efinie de la mani\`ere suivante :

On consid\`ere tout d'abord l'application : $\lambda : {\mathbb C}^{*}
\lra G$ donn\'ee par : 
$\forall t \in {\mathbb C}^{*}, \lambda(t) = \exp \lambda_{t}h^+$, o\`u
$\lambda_{t}$ est un nombre complexe tel que : $e^{\lambda_{t}} =
t$.

Dans ces conditions, on a :
$$
\forall t \in {\mathbb C}^{*}, \forall X \in {\got g}(i), \Ad
\lambda(t).X = t^{i}X
$$
On obtient alors une action de $ {\mathbb C}^{*}$ sur la vari\'et\'e $N$ \`a
l'aide de la formule :
$$
\forall t \in {\mathbb C}^{*}, \forall v \in N, \rho(t). v = t^{-2}\lambda(t).v
$$
On peut, sans pertes de g\'en\'eralit\'e, supposer que chaque $\overline{Z_{i}}$ de la
base consid\'er\'ee est un vecteur de poids $-n_{i}$ et, puisque
le suppl\'ementaire a \'et\'e choisi $\ad h^+$-invariant, que
chaque $X_{j}$ est un vecteur poids, de poids $\nu_{j}$.

On fait agir $ {\mathbb C}^{*}$ sur le coefficient $C_{N}(q)_{l,m}$,
par $\rho$, de la mani\`ere suivante :
$$
\rho(t).C_{N}(q)_{l,m} = C_{N}(\rho(t^{-1}).q)_{l,m}
$$
Ce qui nous donne :
\begin{align*}
\rho(t).C_{N}(q)_{{l,m}} &= C_{N}(\rho(t^{-1}).q)_{{l,m}} \\
&= K(\rho(t^{-1})(e+q), [X_{l},X_{m}]) \\
&= t^{2}K(\lambda(t^{-1})(e+q),[X_{l},X_{m}]) \\
&= t^{2}K(e+q, \lambda(t)[X_{l},X_{m}]) \\
&= t^{2+\nu_{l}+\nu_{m}}K(e+q, [X_{l},X_{m}]) \\
&= t^{2+\nu_{l}+\nu_{m}}C_{N}(q)_{l,m}
\end{align*}
Par ailleurs, suivant la d\'efinition de $\rho$, on a :
$$
\rho(t).(q_{1},\dots,q_{k}) = (t^{-2-n_{1}}q_{1},\dots,t^{-2-n_{k}}q_{k})
$$
Ceci implique que chaque coefficient $C_{N}(q)_{l,m}$ est
quasi-homog\`ene au sens de Slodowy, c'est \`a dire satisfait \`a la
propri\'et\'e suivante :
\begin{equation}
C_{N}(t^{2+n_{1}}q_{1},\dots,t^{2+n_{k}}q_{k})_{l,m} = t^{2+\nu_{l}+\nu_{m}}C_{N}(q_{1},\dots,q_{k})_{l,m}
\end{equation}
 
Soit $\Delta_{N}(q)$ le d\'eterminant de la matrice $C_{N}(q)$. On a :
\begin{align*}
\Delta(\rho(t^{-1}).(q_{1},\dots,q_{k})) &= \sum_{\sigma \in {\mathcal S}_{n}}
\varepsilon_{\sigma} \prod_{i} C_{N}(\rho(t^{-1}).(q_{1},\dots,q_{k}))_{i,\sigma(i)}
\\
&= \sum_{\sigma \in {\mathcal S}_{n}}
\varepsilon_{\sigma} \prod_{i}
t^{2+\nu_{i}+\nu_{\sigma(i)}}C_{N}(q_{1},\dots,q_{k})_{i,\sigma(i)}
\\
&= \sum_{\sigma \in {\mathcal S}_{n}}
\varepsilon_{\sigma} (\prod_{i}
t^{2+\nu_{i}+\nu_{\sigma(i)}}) \prod_{i}C_{N}(q_{1},\dots,q_{k})_{i,\sigma(i)}
\end{align*}
Or, en se servant de la th\'eorie classique associ\'ee \`a la d\'ecomposition
de $ \got g$ en $\got s$-modules, on obtient :
$$
\forall \sigma \in {\mathcal S}_{n}, \sum_{i} \nu_{\sigma(i)} =
\sum_{i} \nu_{i}= - \sum_{k} n_{k} = -p
$$
Il s'en suit que :
$$
\Delta_{N}(t^{2+n_{1}}q_{1}, \dots,t^{2+n_{k}}q_{k}) = \Delta_{N}(q_{1},\dots,q_{k})
$$

Compte-tenu de la d\'efinition des coefficients de $ C_{N}(q)$, le
polyn\^ome $\Delta_{N}(q)$ poss\`ede un terme constant non nul. Ceci
implique donc  que ce d\'eterminant est constant.

Il s'en suit
que la matrice inverse de $C_{N}(q)$ d\'epend polyn\^omialement de la
variable $q$, ce qui d\'emontre le th\'eor\`eme.

\subsection{}
A l'aide de ce qui pr\'ec\`ede, on peut donc d\'efinir, pour chaque
structure transverse $N$ associ\'ee \`a un suppl\'ementaire $ {\got n} \in
{\mathcal N}_{h}$, une notion de degr\'e. Pour cela, on consid\`ere le
degr\'e $d_{N,i,j}(\Lambda)$ du polyn\^ome $\Lambda_{N}(q)_{i,j}$ puis on pose :
$$
deg_{N}(\Lambda) = \sup_{(i,j) \in [1,k]^{2}}deg_{N,i,j}
$$

Dans ([DA]), P.A.Damianou a calcul\'e les degr\'es des structures
transverses des orbites nilpotentes adjointes de $gl(n,{\mathbb C})$,
pour $n \leq 7$, associ\'ees \`a  un certain suppl\'ementaire du
centralisateur.

Il est naturel alors de se poser la question suivante :

{\bf (Q1) :}  L'ensemble $\{deg_{N}(\Lambda), {\got n} \in {\mathcal N}_{h}
\}$ est-il r\'eduit \`a un singleton, c'est \`a dire peut-on dire que le
degr\'e d'une structure transverse est ind\'ependant du choix d'un
suppl\'ementaire $\ad h^+$-invariant ?

La r\'eponse est non comme le montre l' exemple
qui va suivre.

\subsection{}
On suppose que $ {\got g} = sl(4, {\mathbb C})$ et on
consid\`ere l'orbite nilpotente associ\'ee \`a la partition $(3,1)$. On note
$(h^{+},e,f)$ le $sl_{2}$-triplet correspondant comme dans 2.2.  La
caract\'eristique de cette orbite est $(2,0,2)$, sa hauteur est $4$. La
base de racines simples 
$(\alpha_{1}, \alpha_{2}, \alpha_{3})$ v\'erifie alors :
$$
\alpha_{1}(h^{+}) = 2, \ \alpha_{2}(h^{+}) = 0, \ \alpha_{3}(h^{+}) = 2
$$
Dans ces conditions, on peut \'ecrire :
$$
e = X_{\alpha_{1}} + X_{\alpha_{2} + \alpha_{3}}, \ h^{+} =
2H_{\alpha_{1}}+ 2H_{\alpha_{2}}+ 2H_{\alpha_{3}}
$$
Une base de vecteurs-poids de $ {\got g}^{e}$ est donn\'ee par :
\begin{align*}
Z_{1} &= H_{\alpha_{1}} + 2H_{\alpha_{2}} - H_{\alpha_{3}}, & Z_{2} &=
X_{\alpha_{1}} + X_{\alpha_{2} +\alpha_{3}} \\
Z_{3} &= X_{\alpha_{3}}, & Z_{4} &= X_{\alpha_{1}+ \alpha_{2}} \\
Z_{5} &= X_{\alpha_{1}+\alpha_{2}+\alpha_{3}} &
\end{align*}

Consid\'erons ensuite les vecteurs suivants :
\begin{align*}
X_{1} &= H_{\alpha_{1}}, & X_{2} &= H_{\alpha_{2}+ \alpha_{3}}, &
X_{3} &= X_{\alpha_{1}}, &X_{4} = X_{\alpha_{2}} &\\ 
X_{5} &= X_{-\alpha_{2}}, & X_{6} &= X_{-\alpha_{1}}, &
X_{7} &= X_{-\alpha_{3}}, &X_{8} = X_{-\alpha_{1}- \alpha_{2}} &\\
X_{9} &= X_{-\alpha_{2}-\alpha_{3}}, & X_{10} &= X_{-\alpha_{1}-
\alpha_{2}-\alpha_{3}}
\end{align*}
 Soit $ \got n$ le sous-espace vectoriel de $ \got g$ engendr\'e par les
 $(X_{i}, 1 \leq i \leq 10)$. Il s'agit bien d'un suppl\'ementaire $\ad
 h^+$-invariant de $ {\got g}^{e}$. 

La base $(\overline{Z_{i}})$ que nous choisissons pour l'orthogonal $
{\got n}^{\perp}$ est, \`a des constantes non nulles pr\`es, celle qui est issue de
la base de $ \got g$, duale de la base $(Z_{i},X_{j})$ relativement \`a $K$ :
\begin{align*}
\overline{Z_{1}} &= H_{\alpha_{1}} + 2H_{\alpha_{2}} - H_{\alpha_{3}}, & \overline{Z_{2}} &=
X_{-\alpha_{2} -\alpha_{3}} \\
\overline{Z_{3}} &= X_{-\alpha_{3}}, & \overline{Z_{4}} &= X_{-\alpha_{1}- \alpha_{2}} \\
\overline{Z_{5}} &= X_{-\alpha_{1}-\alpha_{2}-\alpha_{3}} &
\end{align*}  

On va  ensuite calculer la matrice $\Lambda'_{N}(q) =
B_{N}(q)C_{N}^{-1}(q)D_{N}(q)$, 
dont le degr\'e nous donnera celui de la
structure transverse $ N = e + {\got n}^{\perp}$, et on obtient :
$$
C_{N}(q) = \begin{pmatrix} 0&0&0&0&0&-2&0&0&1&0 \\
0&0&0&0&0&1&0&0&-2&0 \\0&0&0&q_{4}&0&0&0&0&0&-1 \\
0&0&-q_{4}&0&4q_{1}&0&0&1&0&0 \\0&0&0&-4q_{1}&0&0&-1&0&0&0 \\
2&-1&0&0&0&0&0&0&0&0 \\0&0&0&0&1&0&0&0&0&0 \\
0&0&0&-1&0&0&0&0&0&0 \\-1&2&0&0&0&0&0&0&0&0 \\
0&0&1&0&0&0&0&0&0&0 
\end{pmatrix}
$$

$$
D_{N}(q) = \begin{pmatrix} 0&-q_{2}&0&q_{4}&q_{5} \\
0&2q_{2}&q_{3}&0&q_{5} \\
0&q_{5}&0&0&0 \\
0&-q_{4}&q_{2}&0&0 \\
0&q_{3}&0&0&0 \\
0&0&0&0&q_{2} \\
0&0&4q_{1}&0&-q_{4} \\
0&0&0&-4q_{1}&q_{3} \\
0&0&0&0&0 \\
0&0&0&0&0
\end{pmatrix}
$$

$$
\Lambda'_{N}(q) = \begin{pmatrix} 
0&0&0&0&0 \\
0&0& 4q_{1}q_{3}& -4q_{1}q_{4}& 0 \\
0 & -4q_{1}q_{3}& 0 & 4q_{2}q_{1}- 64 q_{1}^{3}& -
\frac{2}{3}q_{3}q_{2}+ 16q_{1}^{2}q_{3} \\
0&4q_{1}q_{4}&-4q_{1}q_{2}+64q_{1}^{3} & 0 &\frac{2}{3}q_{4}q_{2}-
16q_{1}^{2}q_{4} \\
0&0& \frac{2}{3}q_{3}q_{2} -16 q_{1}^{2}q_{3}&
-\frac{2}{3}q_{4}q_{2} + 16q_{1}^{2}q_{4}&0
\end{pmatrix}
$$
Ceci implique que le degr\'e de la structure transverse $N$ est $3$.

Consid\'erons maintenant le sous-espace
$ {\got n}'$ de $ \got g$ engendr\'e par les
 $(X'_{i}, 1 \leq i \leq 10)$, d\'efinis par :
$$
X'_{1} = H_{\alpha_{1}}, X'_{2} = H_{\alpha_{2}}, X'_{i} = X_{i},
\forall i, 3 \leq i \leq 10
$$
Il s'agit encore  d'un suppl\'ementaire $\ad
 h^+$-invariant de $ {\got g}^{e}$. 

La base $(\overline{Z_{i}})$ de $
{\got n}^{' \perp}$ est alors la suivante :
\begin{align*}
\overline{Z_{1}} &= H_{\alpha_{1}} + 2H_{\alpha_{2}} + 3H_{\alpha_{3}}, & \overline{Z_{2}} &=
X_{-\alpha_{2} -\alpha_{3}} \\
\overline{Z_{3}} &= X_{-\alpha_{3}}, & \overline{Z_{4}} &= X_{-\alpha_{1}- \alpha_{2}} \\
\overline{Z_{5}} &= X_{-\alpha_{1}-\alpha_{2}-\alpha_{3}} &
\end{align*}  
On obtient, cette fois, pour la structure transverse $N' = e + {\got
n}^{'\perp}$ :
$$
C_{N'}(q) = \begin{pmatrix} 0&0&0&0&0&-2&0&0&1&0 \\
0&0&0&0&0&1&0&0&-1&0 \\0&0&0&q_{4}&0&0&0&0&0&-1 \\
0&0&-q_{4}&0&0&0&0&1&0&0 \\0&0&0&0&0&0&-1&0&0&0 \\
2&-1&0&0&0&0&0&0&0&0 \\0&0&0&0&1&0&0&0&0&0 \\
0&0&0&-1&0&0&0&0&0&0 \\-1&1&0&0&0&0&0&0&0&0 \\
0&0&1&0&0&0&0&0&0&0 
\end{pmatrix}
$$

$$
D_{N'}(q) = \begin{pmatrix} 0&-q_{2}&0&q_{4}&q_{5} \\
0&q_{2}&-q_{3}&q_4 &0 \\
0&q_{5}&0&0&0 \\
0&-q_{4}&q_{2}&0&0 \\
0&q_{3}&0&0&0 \\
0&0&0&0&q_{2} \\
0&0&-4q_{1}&0&-q_{4} \\
0&0&0&0&q_{3} \\
0&-4q_1&0&0&0 \\
0&0&0&0&-4q_1
\end{pmatrix}
$$

$$
\Lambda'_{N'}(q) = \begin{pmatrix} 
0&0&0&0&0 \\
0&0& -12q_{1}q_{3}& 12q_{1}q_{4}& 0 \\
0 & 12q_{1}q_{3}& 0 & 0& -2q_{3}q_{2}
 \\
0&-12q_{1}q_{4}&0& 0 &2q_{4}q_{2}  \\
0&0& 2q_{3}q_{2}&
-2q_{4}q_{2}&0
\end{pmatrix}
$$
Ceci nous montre que, dans ce cas, la structure transverse $N'$   est
de degr\'e $2$.

 Il semble, donc, en particulier, que les degr\'es calcul\'es par
Damianou dans ([DA]) ne puissent \^etre consid\'er\'es comme intrins\'eques,
mais plut\^ot d\'ependants du suppl\'ementaire choisi pour les calculs.

\subsection{}

on peut \'egalement se poser la question de savoir si le th\'eor\`eme 2.3
est encore vrai pour un suppl\'ementaire dans $\mathcal N$ qui n'est pas
dans $ {\mathcal N}_{h}$. Si l'on reprend encore une fois l'exemple de
l'orbite sous-r\'eguli\`ere de $sl_{4}({\mathbb C})$, on constate que la
r\'eponse est non.

On consid\`ere donc \`a nouveau l'orbite $(3,1)$ et les notations de
2.5. Soit $ {\got n}_{1}$ le sous-espace vectoriel de $\got g$,
engendr\'e par les vecteurs $(X''_{i}), 1 \leq i \leq 10$, d\'efinis de la
mani\`ere suivante :

\begin{align*}
X''_{1} &= H_{\alpha_{1}}, & X''_{2} &= H_{\alpha_{2}}, &
X''_{3} &= X_{\alpha_{1}}, &X''_{4} = X_{\alpha_{2}} &\\ 
X''_{5} &= X_{-\alpha_{2}}, & X''_{6} &= X_{-\alpha_{1}} + X_{\alpha_{3}}, &
X''_{7} &= X_{-\alpha_{3}}, &X''_{8} = X_{-\alpha_{1}- \alpha_{2}} &\\
X''_{9} &= X_{-\alpha_{2}-\alpha_{3}}, & X''_{10} &= X_{-\alpha_{1}-
\alpha_{2}-\alpha_{3}}
\end{align*}

Il s'agit bien d'un suppl\'ementaire de $ {\got g}^e$ dans $ \got g$;
Cependant, il est facile de v\'erifier que ce suppl\'ementaire n'est pas
$\ad h^{+}$-invariant.

Soit $N_{1} = e + {\got n}_{1}^{\perp}$ la structure transverse
correspondante. La base choisie pour l'orthogonal $ {\got
n}_{1}^{\perp}$ est la suivante :

\begin{align*}
\overline{Z_{1}} &= H_{\alpha_{1}} + 2H_{\alpha_{2}} + 3H_{\alpha_{3}}, & \overline{Z_{2}} &=
X_{-\alpha_{2} -\alpha_{3}} \\
\overline{Z_{3}} &= X_{-\alpha_{3}} - X_{\alpha_{1}}, & \overline{Z_{4}} &= X_{-\alpha_{1}- \alpha_{2}} \\
\overline{Z_{5}} &= X_{-\alpha_{1}-\alpha_{2}-\alpha_{3}} &
\end{align*}

On se place maintenant dans l'ouvert $U_{1} = e + V_{1}$ de $N_{1}$
d\'efini par :
$$
V_{1} = \{ \sum_{i=1}^{i=5}q_{i}\overline{Z_{i}}, q_{3} \not= 1 \}
$$

Les calculs donnent alors :

$$
C_{N_{1}}(q) = \begin{pmatrix} 0&0&0&0&0&-2(1-q_{3})&0&0&1&0 \\
0&0&0&0&0&1-2q_{3}&0&0&-1&0 \\0&0&0&q_{4}&0&0&0&0&0&-1 \\
0&0&-q_{4}&0&0&q_{2}&0&1-q_{3}&0&0 \\0&0&0&0&0&0&-1&0&0&0 \\
2(1-q_{3})&-1+2q_{3}&0&-q_{2}&0&0&4q_{1}&0&0&0 \\0&0&0&0&1&-4q_{1}&0&0&0&0 \\
0&0&0&-1+q_{3}&0&0&0&0&0&0 \\-1&1&0&0&0&0&0&0&0&0 \\
0&0&1&0&0&0&0&0&0&0 
\end{pmatrix}
$$

$$
D_{N_{1}}(q) = \begin{pmatrix} 0&-q_{2}&0&q_{4}&q_{5} \\
0&q_{2}&-q_{3}&q_4 &0 \\
0&q_{5}&0&0&0 \\
0&-q_{4}&q_{2}&0&0 \\
0&q_{3}&0&0&0 \\
4q_{3}&0&0&-q_{5}&q_{2} \\
0&0&-4q_{1}&0&-q_{4} \\
0&0&0&0&q_{3} \\
0&-4q_1&0&0&0 \\
0&q_{3}&0&0&-4q_1
\end{pmatrix}
$$

On v\'erifie que la matrice $C_{N_{1}}(q)$ est inversible, pour $q \in
V_{1}$, et on obtient : 
$$
\Lambda'_{N_{1}}(q) = \begin{pmatrix} 
0&0&4q_{3}^{2}&-8q_{3}q_{4}&-4q_{3}q_{5} \\
0&0& 12q_{1}q_{3}(q_{3}-1)& 12q_{1}q_{4}(1-2q_{3})& -12q_{1}q_{3}q_{5} \\
-4q_{3}^{2} & -12q_{1}q_{3}(q_{3}-1)& 0 & q_{3}q_{5}& -2q_{2}q_{3}
 \\
8q_{1}q_{4}&-12q_{1}q_{4}(1-2q_{3})&-q_{3}q_{5}& 0 & q_{5}^{2} +
2q_{2}q_{4}\frac{2q_{3}-1}{q_{3}-1}  \\
4q_{1}q_{5}&12q_{1}q_{3}q_{5}& 2q_{2}q_{3}&-q_{5}^{2} -
2q_{2}q_{4}\frac{2q_{3}-1}{q_{3}-1}  
&0
\end{pmatrix}
$$
 
On voit bien que cette structure n'est pas polyn\^omiale.

\section{\bf Orbites adjointes conormales : quelques r\'esultats.}

\subsection{}

Supposons, maintenant, qu'il existe, dans $ {\mathcal N}_{h}$ un
suppl\'ementaire $ {\got n}_{0}$ qui soit une alg\`ebre de Lie et
consid\'erons la structure transverse associ\'ee $N_{0} = e + {\got
n}_{0}^{\perp}$.
Suivant la terminologie introduite par M.Roberts dans ([RO]), une
orbite v\'erifiant une telle propri\'et\'e est dite ``{\it conormale}''.   

Il est
clair que, dans ce cas, la matrice $C_{N_{0}}(q)$ ne d\'epend pas de la
variable $q$. Il s'en suit que :
$$
deg_{N_{0}}(\Lambda) \leq 2
$$
On dit alors que la structure transverse $N_{0}$ est ``{\it
quadratique}''. Un tel r\'esultat a \'et\'e donn\'e par Oh dans ([OH]).

Il semble donc int\'eressant de d\'eterminer les orbites
nilpotentes adjointes conormales. On dispose, \`a ce sujet, de deux
r\'esultats que nous allons expliciter dans les paragraphes suivants. 

\subsection{\bf Orbites conormales et orbites sph\'eriques}. 

Consid\'erons la  sous-alg\`ebre de Cartan $ {\got
h}$ de $\got g$, $H$ le sous-groupe de Cartan corres\-pon\-dant dans $G$,
$B$ le 
sous-groupe de Borel de $G$, contenant $H$,  associ\'e \`a un syst\`eme
de racines positives $\Delta^{+}$
choisi dans $\Delta({\got h})$. $U$ sera le radical unipotent de $B$. Soit $\theta$ l'involution de Cartan
sur $G$, telle que l'on ait : $\theta(B) \cap B = H$.

On rappelle alors qu'une orbite
nilpotente est dite ``{\it sph\'erique}'' si cette orbite contient une
$B$-orbite dense. Ces orbites ont \'et\'e classifi\'ees, par exemple par
D.Panyushev dans ([PA 1]) : ce sont celles de hauteur inf\'erieure ou \'egale
\`a $3$. En particulier, les orbites nilpotentes sph\'eriques de
$sl_{n}({\mathbb C})$ sont toutes les orbites de hauteur inf\'erieure
\`a $2$, d\'efinies par les
partitions suivantes : $(2^{a},1^{b})$.

Soit $X$ une $G$-vari\'et\'e irr\'eductible et 
$X^*$ la vari\'et\'e $X$
munie de l'action de $G$ donn\'ee par :
$$\forall g \in G, \forall x^* \in X^*,  g.x^* = \theta(g).x^*$$
$G$ agit alors sur $X \times X^*$ par l'action diagonale, soit :
$g.(x,y) = (g.x,\theta(g).y)$. 

 {\bf D\'efinition : } {\it 

- Un sous-groupe $A$ de $G$ est appel\'e ``stabilisateur en position
g\'en\'erale'' et not\'e s.p.g. s'il existe un ouvert $\Omega$ de $X$ tel que, pour
tout $y \in \Omega, A^y$ est conjugu\'e \`a $A$.

- Dans ce cas, les points de $\Omega$ sont appel\'es ``points en
position g\'en\'erale'' et not\'es p.g.p.}

Du fait de la semi-continuit\'e de la fonction dimension, tout
p.g.p a une $G$-orbite dans $X$ qui est de dimension maximale.

 D.Panyushev, dans ([PA 2], th\'eor\`eme 1.2.2), \'enonce et d\'emontre le r\'esultat
suivant :

 {\bf Th\'eor\`eme 3.1 :} {\it Il existe un point $z = (x,x^*) \in X
\times X^*$ tel que :

1) $U^x$ est un s.p.g pour l'action de $U$ dans $X$.

2) $B^x$ est un s.p.g pour l'action de $B$ dans $X$.

3) $S = G^z = G^x \cap \theta(G^x)$ est un s.p.g pour l'action de
$G$ dans $X \times X^*$.

4)  $U^{x} = U \cap S, \ B^{x} = B \cap S$.

5) Il existe $t \in H, \Delta^+$-dominant, tel que :
\[(G^t)' \subset S \subset G^t\]}
Ici, $(G^t)'$ d\'esigne le sous-groupe d\'eriv\'e de $G^t$.

On suppose maintenant que $X = O$ est une $G$-orbite, on
consid\`ere le point $z = (x,x^*)$ de $O \times O^*$ donn\'e par le
th\'eor\`eme 3.1, appel\'e encore ``point canonique'' et $t
\in H$ l'\'el\'ement satisfaisant au 5). On pose :
$$Q= G^x, L= G^t$$
Dans ce cas, $S = Q \cap \theta(Q)$ est un sous-groupe r\'eductif et
$\theta$-invariant de $Q$.
Notons ${\got b}, {\got l}, {\got q}, {\got s}$ les alg\`ebres de Lie
respectives de $B,L,Q,S$. Soit ${\got t}_0$ l'alg\`ebre de Lie de
$S \cap H$ et ${\got t}_1$ l'orthogonal de ${\got t}_0$ dans $\got t$,
vis-\`a-vis de $K$. D'apr\`es le th\'eor\`eme 3.1, on
sait que : ${\got l} = {\got s} \oplus {\got t}_1$.
On pose  : $P = S.B$.

Comme
$B \cap Q = B \cap S$, on a donc : $P \cap Q = S$.

 {\bf Lemme 3.2 :} {\it $P$ est un sous-groupe parabolique de $G$
contenant $B$.} 

{\bf Preuve :} Posons ${\got
p} = {\got s}+ {\got b}$. On va montrer en fait que, pour toute
racine positive $\alpha$ telle que $X_\alpha \in {\got s}$ et pour
toute racine positive $\beta$ telle que $\beta - \alpha$ soit une
racine, alors le vecteur-racine $X_{\beta - \alpha}$ est dans $\got
p$. Si $\beta-\alpha \geq 0$ , ceci est
clair. Supposons donc $\beta-\alpha \leq 0$. 
On a donc  : $\beta(t) -
\alpha(t) = \beta(t) \leq 0$. Comme $t$ est dominant, ceci implique
que $\beta(t) = 0$ et donc que $X_\beta \in {\got s}$. Ainsi, $ \got
p$ est une alg\`ebre de Lie. Il suffit alors d'utiliser le 5) du
th\'eor\`eme 3.1 pour en d\'eduire que $P$ est un groupe, ce qui d\'emontre le r\'esultat. 

{\bf Th\'eor\`eme  3.3 :} {\it Toute orbite nilpotente sph\'erique est
conormale.}

{\bf Preuve :} On se donne une $G$-orbite nilpotente sph\'erique 
adjointe $O$ et on conserve les notations pr\'ec\'edentes. On peut \'ecrire : ${\got p} =
{\got s} \oplus {\got t}_1 \oplus \ ^u{\got p}$.
D'autre part, on a : $P.x = B.x$ et, comme $x$ est un p.g.p pour
l'action de $B$ dans $O$, l'orbite $B.x$ est de dimension
maximale. Comme, de plus, $O$ est sph\'erique, cette $B$-orbite est
dense dans $G.x$. Il
s'en suit que :
$$
\dim P.x = \dim G.x = \dim {\got g} - \dim {\got q} 
$$
Par ailleurs, on a : 

\begin{align*} 
\dim ({\got q} + {\got p}) &= \dim {\got q} + \dim {\got p}
- \dim {\got s} \\
  &=  \dim {\got q} + \dim {\got p}^x + \dim P.x
- \dim {\got s} \\
&= \dim {\got g} + \dim{\got p}^{x} - \dim {\got s}
\end{align*}

Comme ${\got s} = {\got p}^x$, on en d\'eduit que :
$$
\dim ({\got q} + {\got p}) = \dim {\got g}
$$
On a, donc, pour une telle
orbite :\[{\got g} = {\got q} + {\got p} = {\got q} \oplus ({\got t}_1
\oplus \ ^u{\got p})\]
Posons : ${\got s}_O = {\got t}_1 \oplus \ ^u{\got p}$. 
Alors, il est clair que ${\got s}_O$ est une sous-alg\`ebre de $\got
g$ et que $\got q$ est le centralisateur dans $\got g$ d'un point de
$O$, ce qui d\'emontre le th\'eor\`eme 3.3. 

\subsection{\bf Orbites conormales dans $\bf sl_{n}({\mathbb C})$.}
Nous nous int\'eressons maintenant aux orbites conormales de
$sl_{n}({\mathbb C})$.  Certaines de ces orbites
ont d\'ej\`a \'et\'e exhib\'ees par M.Rais dans ([RA]). Le m\^eme auteur a
aussi prouv\'e, dans des notes non publi\'ees, que toute orbite nilpotente
minimale et toute orbite nilpotente r\'eguli\`ere  d'une alg\`ebre de Lie
semi-simple quelconque  \'etaient  conormales.

Nous nous proposons de d\'emontrer le r\'esultat suivant :

 {\bf Th\'eor\`eme 3.4 :} {\it Soit $O$ une orbite nilpotente de
$sl_{n}({\mathbb C})$ de partition $(p_{1},p_{2}, \dots,p_{s})$ telle
que :
$$
\forall i,j, 1 \leq i,j \leq s, \ |p_{i}-p_{j}| \leq 1
$$
Alors, $O$ est conormale.}

 Il faut noter que cette classe d'orbites nilpotentes contient les
orbites sph\'eriques et l'orbite r\'eguli\`ere.

La d\'emonstration de ce th\'eor\`eme est essentiellement technique et
repose sur la d\'e\-termi\-na\-tion explicite du centralisateur d'un \'el\'ement
convenablement choisi de l'orbite.

On reprend \`a ce sujet les notations introduites dans 2.2. La
d\'etermination des caract\'eristiques de la famille d'orbites
nilpotentes  de $sl_n$ donn\'ee dans l'\'enonc\'e du
th\'eor\`eme 3.4 permet de d\'ecomposer celle-ci selon  les deux types de partitions
suivants :
\begin{itemize}
\item Type I. $(p^{r}), \ pr = n$. La caract\'eristique est alors :
$(B^{p-1},0^{r-1})$ o\`u  $B =(0^{r-1},2)$.
\item Type II. $(p^r,(p-1)^s),  \  rp + s(p-1)= n$. La caract\'eristique
est alors : $(B^{p-1},0^{r-1})$ o\`u  $B = (0^{r-1},1,0^{s-1},1)$. 
\end{itemize}

\subsection{Preuve du th\'eor\`eme 3.4 : le cas du type I.} Soit $O = (p^{r})$ une
orbite de type I.
L'ensemble $\Pi({\got h})$ des racines simples peut se d\'ecrire de la mani\`ere
suivante :
\begin{align*}
\alpha_{kr}(h^{+}) &= 2, \forall k, \ 1 \leq k \leq p-1 \\
\alpha_{kr+i}(h^{+}) &= 0, \forall (k,i) \in \  0 \leq k \leq
p-1, \  1\leq i \leq 
r-1 \\
\Pi({\got h})& = \{\alpha_{kr}, \forall k, \ 1\leq k \leq 
p-1\} \cup \{\alpha_{kr+i}, \forall (k,i) \ 0\leq k \leq 
p-1, \  1 \leq i \leq 
r-1\}
\end{align*}

Posons :
\begin{align*}
\beta_{i,j} &= \sum_{s=i}^{s=j} \alpha_s , \forall i,j, \  i \leq j\\
X(i,j) &= X_{\beta_{i,j}} 
\end{align*}
D'autre part, nous noterons dor\'enavant $sl_q(\Pi')$ la
sous-alg\`ebre de $\got g$, isomorphe \`a $sl_q({\mathbb C})$, engendr\'ee par
le syst\`eme de racines simples $\Pi', \Pi' \subset \Pi({\got h}), \ \natural \Pi' = q-1$.

Consid\'erons maintenant la d\'ecomposition de $\got g$ en
sous-espaces propres suivant l'action de $\ad h^{+}$, soit :
\[ {\got g} = \bigoplus_{m=-p + 1}^{m= p- 1} {\got g}(2m)\] 
le calcul nous donne, pour tout entier $m$ tel que $1 \leq m \leq p-1$ :
\begin{align*}
{\got g}(0) &= \bigoplus_{k=0}^{k=p-1} sl_p(\alpha_{kr+i}, 1 \leq
i \leq r-1) \oplus < H_{\alpha_{kr}}, 1 \leq k \leq p-1 > \\
{\got g}(2m) &= < X(kr+i,(k+m)r+j),  0 \leq k \leq p-m-1, \ 1 \leq
i \leq r, \ 0 \leq j \leq r-1 > 
\end{align*}

Le g\'enerateur $e$ de l'orbite est alors donn\'e par :
$$
e = \sum_{i=1}^{i= (p-1)r} X(i,p-1+i)
$$

Venons-en au stabilisateur ${\got g}^e$ de $e$ dans ${\got g}$. On a
la d\'ecomposition :
$$
{\got g}^e = \bigoplus_{m=0}^{m=p-1}{\got g}^e(2m)
$$
Consid\'erons les vecteurs suivants :
\begin{align*}
\forall i, 1 \leq i \leq r-1, X^{+}_{i} &=
\sum_{k=0}^{k=p-1}X_{\alpha_{kr+i}}\\  
X^{-}_{i} &=
\sum_{k=0}^{k=p-1}X_{-\alpha_{kr+i}}\\
H_{i} &=
\sum_{k=0}^{k=p-1}H_{\alpha_{kr+i}}
\end{align*}
La famille ${\mathcal F}_{0}(e) = (X^{+}_{i},X^{-}_{i},H_{i}, 1 \leq i \leq r-1)$ engendre
une sous-alg\`ebre  de $\got g$, isomorphe \`a $sl_{r}({\mathbb
C})$. Nous la noterons $sl_{r}({\mathcal F}_{0}(e))$.

Le calcul nous donne alors, pour tout entier $m$ tel que $1 \leq m \leq p-1$ :
\begin{align*}
{\got g}^e(0) &= sl_r({\mathcal F}_{0}(e)) \\
{\got g}^e(2m) &= < \sum_{k=0}^{k=p-1}X(kr+i,(k+m)r+j), \ 1 \leq
i \leq r, \ 0 \leq j \leq r-1 > 
\end{align*} 

L'ensemble $\{ \alpha_{kr+i}, \  0 \leq k \leq p - 2, 1 \leq i
\leq r-1\} \cup \{\alpha_{kr}, 1 \leq k \leq p-2\}$ d\'efinit un sous-syst\`eme de racines de $\Pi({\got h})$, que nous
noterons $\Pi_p({\got h})$. 

 Posons :
$$
{\got g}_p = sl_{(p-1)r}(\Pi_p({\got h}))
$$
Consid\'erons maintenant la sous-alg\`ebre parabolique maximale ${\got
p}_p$ de ${\got
g}$, obtenue \`a partir du syst\`eme de racines simples
$\Pi \backslash \{\alpha_{(p-1)r}\}$. Soit $ {\got p}^{-}_{p}$ la
sous-alg\`ebre parabolique oppos\'ee et $^{u}{\got p}^{-}_{p}$ son radical
unipotent.
On a : ${\got
p}_p^{-} = {\got m}_p \oplus {\got a}_p \oplus \ ^u{\got
p}_p^{-}$, avec :
\begin{itemize}
\item ${\got m}_p = {\got g}_p \oplus
sl_p(\alpha_{(p-1)r + i}, 1 \leq i \leq r-1)$ \\
\item ${\got a}_p = <H_{\alpha_{(p-1)r}}>$
\end{itemize}
Soit, enfin : 
$$
{\got s}_p = {\got g}_p \oplus {\got
a}_p \oplus \ ^u{\got p}_p^{-} 
$$
Il est clair, alors,  que ${\got s}_p$ est une sous-alg\`ebre de $\got g$, que
cette alg\`ebre est $\ad h^{+}$-invariante 
et que : ${\got g} = {\got g}^e \oplus {\got s}_p$. 

Ceci nous donne le r\'esultat souhait\'e, en ce qui concerne le type I.

\subsection{Preuve du th\'eor\`eme 3.4 : le cas du type II}

Soit $O = (p^r,(p-1)^s)$
une orbite de type II. 
Posons : $T = r+s$.

 L'ensemble $\Pi({\got h})$ des racines simples se d\'ecrit alors de
la mani\`ere suivante :
\begin{align*}
\alpha_{kT + r}(h^{+}) &= 1, \forall k, \ 0\leq k \leq p-2  \\
\alpha_{(k+1)T}(h^{+}) &= 1, \forall k, \ 0\leq k \leq p-2 
\\
\alpha_{kT + i}(h^{+}) &= 0, \forall (k,i), \ 0\leq k \leq p-1,
1\leq i \leq r-1 \\
\alpha_{kT + i}(h^{+}) &= 0, \forall (k,i), \ 0\leq k \leq p-2, \ 
r+1\leq i \leq T-1 \\
\Pi({\got h}) &= \{\alpha_{kT + i}, (k,i), \ 0\leq k \leq p-2,
1\leq i \leq T \} \cup \{ \alpha_{(p-1)T
+ i}, \ 1\leq i \leq r-1\}
\end{align*}  

On consid\`ere ensuite la d\'ecomposition de $\got g$ en sous-espaces
propres suivant l'action de $ \ad h^{+}$, soit :
$$
 {\got g} = \bigoplus_{m=-p + 1}^{m= p - 1} {\got g}(2m)
\oplus {\got g}(2m-1) 
$$
le calcul nous donne, pour tout entier $m, 1 \leq m \leq p-1$ :

\begin{align*}
{\got g}(0) = < H_{\alpha_{kT+r}}, H_{\alpha_{(k+1)T}}, 0 \leq k \leq
p-2 > \oplus  &\bigoplus_{k=0}^{k=p-1} sl_{r}(\alpha_{kT+i}, 1 \leq
i \leq r-1) \\ 
&\bigoplus_{k=0}^{k=p-2} sl_{s}(\alpha_{kT+i}, r+1 \leq i
\leq T-1) 
\end{align*} 
\begin{alignat*}{2} 
{\got g}(2m) &= < X(kT+i,(k+m)T+j),&  &\ 1 \leq i \leq r, 0 \leq j \leq
r-1, \\
& & &0 \leq k \leq p-m-1 >\\ 
&\oplus  < X(kT+i,(k+m)T+j),& &\  r+1 \leq i \leq T,
r \leq j \leq T-1, \\ 
& & &0 \leq k \leq p-m-1 >
\end{alignat*}
\begin{alignat*}{2}
{\got g}(2m-1) &= < X(kT+i,(k+m-1)T+j),& &  \  1 \leq
i  \leq r, r \leq j \leq T-1, \\
& & &0 \leq k \leq p-m-1> \\
 &\oplus < X(kT+i,(k+m)T+j), & & r+1 \leq i  \leq T, 0
\leq j \leq r-1, \\
& & &0 \leq k \leq p-m-1> 
\end{alignat*}

Le g\'en\'erateur $e$ de l'orbite est alors donn\'e par :
$$
e = \sum_{i=1}^{i=n-T} X(i,i+T-1)
$$

Consid\'erons les vecteurs suivants :
\begin{align*}
X^{+}_{i} &=
\sum_{k=0}^{k=p-1}X_{\alpha_{kT+i}}\\  
X^{-}_{i} &=
\sum_{k=0}^{k=p-1}X_{-\alpha_{kT+i}}\\
H_{i} &=
\sum_{k=0}^{k=p-1}H_{\alpha_{kT+i}}
\end{align*}
Les familles ${\mathcal F}_{0,r}(e) = (X^{+}_{i},X^{-}_{i},H_{i}, 1
\leq i \leq r-1)$ et ${\mathcal F}_{0,s}(e) = (X^{+}_{i},X^{-}_{i},H_{i}, r+1
\leq i \leq T-1)$ engendrent des 
sous-alg\`ebres de $\got g$, isomorphes respectivement \`a $sl_{r}({\mathbb
C})$ et $sl_{s}({\mathbb C}) $. Nous les  noterons $sl_{r}({\mathcal
F}_{0,r}(e))$ et $sl_{s}({\mathcal F}_{0,s}(e))$.

Venons-en au stabilisateur ${\got g}^e$ de $e$ dans ${\got g}$. On a
la d\'ecomposition :
$$
 {\got g}^e = \bigoplus_{m=0}^{m=p-1}{\got g}^e(2m) \oplus {\got g}^e(2m-1)
$$
On notera $ {\got z}_{0}(e)$ le centre de $ {\got g}^e(0)$. Pour des
raisons de dimension, on v\'erifie qu'il s'agit d'un sous-espace de
dimension $1$ de $ \got h$. 

Le calcul nous donne, pour tout entier $m, 1 \leq m \leq p-1$ :
\begin{align*}
{\got g}^e(0) &= sl_{r}({\mathcal F}_{0,r}(e)) \oplus
sl_{s}({\mathcal F}_{0,s}(e)) \oplus {\got z}_{0}(e)  \\
{\got g}^e(2m) &= < \sum_{k=0}^{k=p-m-1} X(kT+i,(k+m)T+j), 
 1 \leq
i \leq r, 0 \leq j \leq r-1> \\
 & \oplus < \sum_{k=0}^{k=p-m-1} X(kT+i,(k+m)T+j), r+1 \leq i \leq T, r
\leq j \leq T-1 >\\
{\got g}^e(2m-1) &= < \sum_{k=0}^{k=p-m-1}X(kT+i,(k+m-1)T+j),  1\leq
i \leq r, r \leq j \leq T-1 >\\
& \oplus < \sum_{k=0}^{k=p-m-1}X(kT+i,(k+m)T+j), r+1 \leq i \leq T, 0 \leq j \leq
r-1> \\
\end{align*}

L'ensemble $\{ \alpha_{kT+i}, \  0 \leq k \leq p - 2, 1 \leq i
\leq r-1\} \cup \{ \alpha_{kT+i}, \  0 \leq k \leq p - 3, r \leq i
\leq T\}$  d\'efinit un sous-syst\`eme de racines simples  de
$\Pi({\got h})$, de cardinal $n-T-1$, que nous
noterons $\Pi_p({\got h})$. Posons, suivant les notations introduites
pr\'ecedemment :
\[{\got g}_p = sl_{n-T}(\Pi_p)\]

Consid\'erons maintenant la sous-alg\`ebre parabolique ${\got
p}_p$ de ${\got
g}$, obtenue \`a partir du syst\`eme de racines simples
$\Pi \backslash \{\alpha_{(p-2)T + r},\alpha_{(p-1)T}\}$ et $ {\got
p}_{p}^{-}$ la sous-alg\`ebre parabolique oppos\'ee. On a : ${\got
p}_p^{-} = {\got m}_p \oplus {\got a}_p \oplus \ ^u{\got
p}_p^{-}$, avec :
\begin{align*}
{\got m}_p &= {\got g}_p \oplus
sl_s(\alpha_{(p-2)T + i}, r+1 \leq i \leq T-1) 
 \oplus sl_r(\alpha_{(p-1)T + i}, 1 \leq i \leq r-1) \\
{\got a}_p &= <H_{\alpha_{(p-1)T}}, H_{\alpha_{(p-2)T
+ r}}>
\end{align*}

On choisit, ensuite, un suppl\'ementaire quelconque de $ {\got g}^e
\cap {\got a}_{p}$ dans $ {\got a}_{p}$. Il s'agit d'un espace de
dimension $1$, que nous noterons $\widetilde{{\got a}_{p}}$.

Soit, enfin : ${\got s}_p = {\got g}_p \oplus {\widetilde{{\got
a}_p}} \oplus \ ^u{\got p}_p$. 

On v\'erifie  que ${\got s}_p$ est une sous-alg\`ebre de $\got g$, que
cette alg\`ebre est $\ad h^{+}$-invariante 
et que : ${\got g} = {\got g}^e \oplus {\got s}_p$.

Ceci d\'emontre compl\`etement le th\'eor\`eme 3.4.

Pour finir, Je formule \'egalement la conjecture suivante :

{\bf Conjecture :} {\it Les seules orbites nilpotentes conormales
de $sl_{n}({\mathbb C})$ sont celles donn\'ees par le th\'eor\`eme 3.3.}

\subsection{\bf Retour aux structures transverses quadratiques.}

L'exemple de l'orbite $(3,1)$ du paragraphe 2.5 nous montre
qu'il existe des structures transverses quadratiques \`a une orbite
adjointe correspondant \`a un suppl\'ementaire qui n'est pas une alg\`ebre
de Lie.

On peut alors se poser la question suivante :

{\bf (Q2) :} La structure transverse \`a  une orbite nilpotente
adjointe conormale est-elle toujours quadratique ? 

Encore une fois, la r\'eponse est non comme le montre l'exemple suivant :

On consid\`ere cette fois l'orbite nilpotente de partition $(3,2)$
dans $ {\got g} = sl_{5}({\mathbb C})$. Cette orbite est conormale,
d'apr\`es le th\'eor\`eme 3.3. Soit $(h^+,e,f)$ le $sl_{2}$-triplet associ\'e,
$(1,1,1,1)$ la caract\'eristique, $(\alpha_{i}), 1 \leq i \leq 4$, le
syst\`eme de racines simples correspondant. On a :
$$
e = X_{\alpha_{1}+ \alpha_{2}} +  X_{\alpha_{2}+ \alpha_{3}} +
X_{\alpha_{3}+ \alpha_{4}} 
$$
Une base de vecteurs-poids de $ {\got g}^{e}$ est donn\'ee par :

\begin{align*}
Z_{1} &=2 H_{\alpha_{1}} -H_{\alpha_{2}} + H_{\alpha_{3}} -
2H_{\alpha_{4}}, & Z_{2} &= 
X_{\alpha_{1}} + X_{\alpha_{3}} &\\  
Z_{3} &= X_{\alpha_{2}} + X_{\alpha_{4}}, &
Z_{4} &= X_{\alpha_{1}+
\alpha_{2}}+ X_{\alpha_{3} + \alpha_{4}} & \\
Z_{5} &= X_{\alpha_{2} +
\alpha_{3}}, &
Z_{6} &= X_{\alpha_{1}+\alpha_{2}+\alpha_{3}} &\\
Z_{7} &= X_{\alpha_{2}+\alpha_{3}+\alpha_{4}}, & Z_{8}&=
X_{\alpha_{1}+ \alpha_{2}+ \alpha_{3}+\alpha_{4}}   
\end{align*}

Soit $ {\got n}_{0}$ la sous-alg\`ebre de $ \got g$ donn\'ee par la
d\'emonstration du th\'eor\`eme 3.3, telle que : ${\got
g} = {\got g}^{e} \oplus {\got n}_{0}$. La structure transverse $N_{0}
= e + {\got n}_{0}^{\perp}$ est donc quadratique.

Consid\'erons maintenant le suppl\'ementaire $ {\got n}_{f} = Im \ad f$ de $ {\got
g}^{e}, \ {\got n}_{f}^{\perp} = {\got g}^{f}$ et $N_{f} = e + {\got g}^{f}$. On utilise la base
suivante de $ {\got n}_{f}$ :
\begin{align*}
X_{1} &= X_{-\alpha_{1}-\alpha_{2}-\alpha_{3}-\alpha_{4}}, & X_{2} &=
X_{-\alpha_{1}-\alpha_{2}-\alpha_{3}},& X_{3} &=
X_{-\alpha_{2}-\alpha_{3}-\alpha_{4}}, & X_{4} &=
X_{-\alpha_{1}-\alpha_{2}} \\
 X_{5} &= X_{-\alpha_{2}-\alpha_{3}}, &
X_{6} &= X_{-\alpha_{3}-\alpha_{4}}, & X_{7} &= X_{-\alpha_{1}}, &
X_{8} &= X_{-\alpha_{2}} \\
X_{9} &= X_{-\alpha_{3}}, & X_{10} &= X_{-\alpha_{4}}, & X_{11} &=
H_{\alpha_{1}+\alpha_{2}} \\
 X_{12} &= H_{\alpha_{2}+\alpha_{3}},
&X_{13} &= H_{\alpha_{3}+\alpha_{4}}, & X_{14} &=
X_{\alpha_{1}}-X_{\alpha_{3}} \\
X_{15} &=
X_{\alpha_{2}}-X_{\alpha_{4}}, & X_{16} &=
X_{\alpha_{1}+ \alpha_{2}}-X_{\alpha_{3}+ \alpha_{4}}
\end{align*}

La base choisie, pour $ {\got g}^{f}$, est la suivante :
\begin{align*}
\overline{Z_{1}} &=2 H_{\alpha_{1}} -H_{\alpha_{2}} + H_{\alpha_{3}} -
2H_{\alpha_{4}}, & \overline{Z_{2}} &= 
X_{-\alpha_{1}} + X_{-\alpha_{3}} &\\  
\overline{Z_{3}} &= X_{-\alpha_{2}} + X_{-\alpha_{4}}, &
\overline{Z_{4}} &= X_{-\alpha_{1}
-\alpha_{2}}+ X_{-\alpha_{3} - \alpha_{4}} & \\
\overline{Z_{5}} &= X_{-\alpha_{2} 
-\alpha_{3}}, &
\overline{Z_{6}} &= X_{-\alpha_{1}-\alpha_{2}-\alpha_{3}} &\\
\overline{Z_{7}} &= X_{-\alpha_{2}-\alpha_{3}-\alpha_{4}}, & \overline{Z_{8}}&=
X_{-\alpha_{1}- \alpha_{2}- \alpha_{3}-\alpha_{4}}   
\end{align*}

 On reprend ensuite les notations pr\'ec\'edentes et le calcul nous donne :      

$$
\Lambda_{N}'(q) = \begin{pmatrix}
0&0&0&0&0&0&0&0 \\
0&0& -75q_{1}^{2}& 15q_{1}q_{2}& -5q_{1}q_{2}&
\frac{1}{3}q_{2}^{2}& \{z_{2},z_{7}\}&\{z_{2},z_{8}\} \\
0&75q_{1}^{2}&0&-15q_{3}q_{1}&5q_{3}q_{1}&
\{z_{3},z_{6}\}&-\frac{1}{3}q_{3}^{2}&\{z_{3},z_{8}\} \\

0&-15q_{1}q_{2}&15q_{3}q_{1}&0&0&
\{z_{4},z_{6}\} &\{z_{4},z_{7}\} &\{z_{4},z_{8}\} \\
0&5q_{1}q_{2}&-5q_{3}q_{1}&0&0&\{z_{5},z_{6}\}&\{z_{5},z_{7}\}&
\{z_{5},z_{8}\} \\
0&-\frac{1}{3}q_{2}^{2}&-\{z_{3},z_{6}\}&-\{z_{4},z_{6}\}&-\{z_{5},z_{6}\}&
0&\{z_{6},z_{7}\}&\{z_{6},z_{8}\} \\
0& - \{z_{2},z_{7}\}& - \{z_{3},z_{7}\}&- \{z_{4},z_{7}\}&-
\{z_{5},z_{7}\}&
- \{z_{6},z_{7}\}&0& \{z_{7},z_{8}\} \\
0& - \{z_{2},z_{8}\}& - \{z_{3},z_{8}\}&- \{z_{4},z_{8}\}&-
\{z_{5},z_{8}\}&
- \{z_{6},z_{8}\}& -\{z_{7},z_{8}\}&0  
\end{pmatrix}
$$
\begin{align*}
\{z_{2},z_{7}\} &= \frac{375}{2}q_{1}^{3} + 10q_{5}q_{1}-5q_{4}q_{1}
- \frac{2}{3}q_{2}q_{3} \\ 
\{z_{2},z_{8}\} &= -\frac{75}{2}q_{1}^{2}q_{2} + 10q_{1}
q_{6} -\frac{1}{3}q_{2}q_{4} \\
\{z_{3},z_{6}\} &= \frac{125}{2}q_{1}^{3} - 10q_{5}q_{1}+5q_{4}q_{1}
+ \frac{2}{3}q_{2}q_{3} \\  
\{z_{3},z_{8}\} &= -\frac{25}{2}q_{1}^{2}q_{3}
-10q_{1}q_{7}+ \frac{1}{3}q_{3}q_{4} \\
\{z_{4},z_{6}\} &=
-\frac{25}{2}q_{1}^{2}q_{2} + 2q_{2}q_{5}-\frac{5}{2}q_{1}q_{6} 
\end{align*}
\begin{align*}
\{z_{4},z_{7}\}&= -\frac{75}{2}q_{1}^{2}q_{3} -2q_{3}q_{5} +
\frac{5}{2}q_{1}q_{7} \\
\{z_{4},z_{8}\} &= -\frac{3}{2}q_{3}q_{6}+\frac{3}{2}q_{2}q_{7}+
 10q_{1}q_{2}q_{3} \\
\{z_{5},z_{6}\} &=
\frac{25}{2}q_{1}^{2}q_{2} - 2q_{2}q_{5}- \frac{5}{2}q_{1}q_{6} \\ 
\{z_{5},z_{7}\}&= \frac{25}{2}q_{1}^{2}q_{3} +2q_{3}q_{5} +
\frac{5}{2}q_{1}q_{7} \\ 
\{z_{5},z_{8}\} &=\frac{3}{2}q_{3}q_{6}-\frac{3}{2}q_{2}q_{7}-
 5q_{1}q_{2}q_{3} \\
\{z_{6},z_{7}\}&= \frac{625}{4}q_{1}^{4} 
-\frac{25}{2}q_{1}^{2}q_{4} -
5q_{1}q_{2}q_{3} -25q_{1}^2q_5 + 5q_{1}q_{8}  
+ \frac{1}{6}q_{2}q_{7}  + 2q_{4}q_{5} - q_{5}^{2} +
\frac{1}{6}q_{3}q_{6} \\
\{z_{6},z_{8}\} &=-\frac{125}{4}q_{1}^{3}q_{2} +
\frac{5}{2}q_{1}q_{2}q_{4} + 5q_1q_2q_5  
+ \frac{3}{2}q_{2}^{2}q_{3} -
\frac{3}{2}q_{2}q_{8} - q_{5}q_{6} + \frac{4}{3}q_{4}q_{6}\\
\{z_{7},z_{8}\} &=\frac{125}{4}q_{1}^{3}q_{3}
-\frac{15}{2}q_{1}q_{3}q_{4} + 25q_1^2q_7 
- \frac{3}{2}q_{2}q_{3}^{2}+
\frac{3}{2}q_{3}q_{8} - \frac{4}{3}q_{4}q_{7}+ q_{5}q_{7}
\end{align*}
 On constate, alors, que le degr\'e de la structure transverse $N$ est $4$.

{\bf R\'ef\'erences.}

{\bf [BE-GO]) : P.G.Bergmann,I.Goldberg} {\it Dirac bracket
transformations in phase space.} Phys.Rev. {\bf 98} (1955) 531-538

{\bf [CU-RO] : R.Cushman,M.Roberts} {\it Poisson structures
transverse to coadjoint orbits.} Bull.Sci.math {\bf 126} (2002) 525-534.

{\bf [DA] : P.A.Damianou.}{\it Transverse Poisson structures of
coadjoint orbits.} Bull. Sci. math. {\bf 120} (1996),  195-214.

{\bf [DI] : J.Dixmier} {\it Enveloping algebras.} Graduate Studies
in Mathematics, Vol.{\bf 11}.

{\bf [OH] : G.Y.Oh.} {\it Some remarks on the transverse Poisson
structures of coadjoint orbits.} Lett.Math.Phys. {\bf 12} (1986), 87-91

{\bf [PA 1] : D.Panyushev.} {\it On spherical nilpotent orbits and
beyond}.
Ann.Inst.Fourier, {\bf 49}, 5, (1999), 1453-1476

{\bf [PA 2] : D.Panyushev.} {\it Complexity and rank of actions in
invariant theory.} Journ. of Math. Science, 95, 1 (1999) 1925-1985

{\bf [RA] : M.Rais.} {\it La repr\'esentation adjointe du groupe
affine}. Ann.Inst.Fourier,{\bf 28}, 1 (1978), 207-237.  

{\bf [RO] : M.Roberts, C.Wulff, JSW Lamb.} {\it Hamiltonian
systems near relative equilibria.} J.Differentiel equations. {\bf 179}
(2002), 562-604.
 
{\bf [SG] : M.Saint-Germain} {\it Poisson algebras and transverse
structures.} J.Geom.Phys. {\bf 31} (1999), 153-194.

{\bf [SL] : P.Slodowy.} {\it Simple singularities and Simple
algebraic groups.} Lecture Notes in Math. {\bf 815}. Springer Verlag
Berlin, Heidelberg, New-York. 1980.

{\bf [SP-ST] : T.A.Springer,R.Steinberg} {\it Conjugacy classes.} in :
``seminar on algebraic groups and related finite groups'',
ed. A.Borel, et.al., 168-266, Lecture notes in Math. {\bf 131}
Springer-Verlag Berlin, Heidelberg, New-York. 1970

{\bf [WE] : A.Weinstein.} {\it Local structure of Poisson
manifolds.} J.Differential Geom. {\bf 18} (1986), 523-557

\end{document}